\documentclass[11pt, reqno]{amsart}
\usepackage{amsmath, amssymb, amsthm, verbatim,enumerate,bbm, mathtools,color, enumitem}
\usepackage{dsfont}
\usepackage{hyperref}
\usepackage{lineno}
\usepackage{fullpage}

\newtheorem{theorem}{{Theorem}}

\newtheorem{rem}[theorem]{{Remark}}

\newtheorem{thm}[theorem]{{Theorem}}

\newtheorem{lem}[theorem]{{Lemma}}
\newtheorem{defi}[theorem]{{Definition}}

\numberwithin{equation}{section}
\numberwithin{theorem}{section}

\newcommand{\var}{\varepsilon}

\title{Dirac-type Problem of Rainbow matchings and Hamilton cycles in Random Graphs}
\author{Asaf Ferber}
\address{Department of Mathematics, University of California, Irvine.
Email:  {{asaff@uci.edu}}. }
\thanks{The first author is supported in part by NSF grant DMS-1953799, NSF Career DMS-2146406, and a Sloan's fellowship. The second author is partially supported by National Natural Science Foundation of China (12371341).} 
\author{Jie Han}
\address{School of Mathematics and Statistics and Center for Applied Mathematics, Beijing Institute of Technology, Beijing, China. Email: {{han.jie@bit.edu.cn}}} 
\author{Dingjia Mao}
\address{Department of Mathematics, University of California, Irvine.
Email:  {{dingjiam@uci.edu}}}
\date{\today}

\begin{document}
\maketitle

\begin{abstract}


Given a family of graphs $G_1,\dots,G_{n}$ on the same vertex set $[n]$, a rainbow Hamilton cycle is a Hamilton cycle on $[n]$ such that each $G_c$ contributes exactly one edge. We prove that if $G_1,\dots,G_{n}$ are independent samples of $G(n,p)$ on the same vertex set $[n]$, then for each $\var>0$, whp, every collection of spanning subgraphs $H_c\subseteq G_c$, with $\delta(H_c)\geq(\frac{1}{2}+\var)np$, admits a rainbow Hamilton cycle.
A similar result is proved for rainbow perfect matchings in a family of $n/2$ graphs on the same vertex set $[n]$.

\end{abstract}


\section{Introduction}\label{s:intro}

\subsection{Dirac-type problems}
Arguably the two most studied objects in graph theory are \emph{perfect matchings} and \emph{Hamilton cycles}. A perfect matching in a graph $G=(V,E)$ is a collection of vertex-disjoint edges which covers $V$, and a Hamilton cycle is a cycle passing through all the vertices of $G$. For ensuring a perfect matching, the celebrated Hall's theorem \cite{hall1935representatives} and Tutte's theorem \cite{tutte1956theorem} each provides an equivalent statement. However, as opposed to the problem of finding a perfect matching (if one exists) in a graph $G$ which has efficient (polynomial-time) resolutions, the analogous problem for Hamilton cycles is listed as one of the NP-hard problems by Karp \cite{karp1972reducibility}.  Therefore, as one cannot hope to find a Hamilton cycle efficiently, it is natural to study sufficient conditions which guarantee its existence.

One of the first results of this type is the celebrated theorem by Dirac \cite{dirac1952some}, which states that every graph on $n\geq 3$ vertices with minimum degree $n/2$ is \emph{Hamiltonian}, that is, contains a Hamilton cycle (and in particular, if $n$ is even, then it also contains a perfect matching). While Dirac's theorem is sharp in general, one would like to find sufficient conditions for sparser graphs. A natural candidate to begin with is a typical graph sampled from the binomial random graph model $G(n,p)$. That is, a graph $G$ on vertex set $[n]$, where each (unordered) pair is being sampled as an edge with probability $p$, independently. In 1960, Erd\H{o}s and R\'enyi raised a question of what the threshold probability of Hamiltonicity in random graphs is. This question attracted a lot of attention in the past few decades. After a series of efforts by various researchers, including Korshunov \cite{korshunov1976solution} and P\'osa \cite{posa1976hamiltonian}, the problem was finally solved
by Koml\'os and Szemer\'edi \cite{komlos1983limit} and independently by Bollob\'as \cite{bollobas1984evolution}, who proved that if
$p \geq (\log_{} n + \log_{} \log_{} n + \omega(1))/n$, where $\omega(1)$ tends to infinity with $n$ arbitrarily slowly,
then the probability of the random graph $G(n, p)$ being Hamiltonian tends to~1 (we say such an event happens \emph{with high probability}, or \emph{whp} for brevity). This result is best possible, since for $p \leq (\log_{} n + \log_{} \log_{} n - \omega(1))/n$ whp there are vertices of
degree at most one in $G(n, p)$ (see, e.g. \cite{bollobas1998random}). More precisely, they showed that for $p=\frac{\log n+\log\log n+c_n}{n}$, 
\[
\lim_{n\to\infty} \Pr[G(n,p)\text{ has a Hamilton cycle}]=
\begin{cases}
    0&\text{ if }c_n\to-\infty\\
    e^{-e^{-c}}&\text{ if } c_n\to c\\
    1&\text{ if }c_n\to\infty
\end{cases}.
\] 
Not surprisingly, the probability that the random graph $G(n,p)$ contains a perfect matching is also precise asymptotically. In fact, Erd\H os and R\'enyi \cite{erdHos1966existence} proved that for $p=\frac{\log n+c_n}{n}$, 
\[
\lim_{n\to\infty,\: n\text{ even}} \Pr[G(n,p)\text{ has a perfect matching}]=
\begin{cases}
    0&\text{ if }c_n\to-\infty\\
    e^{-e^{-c}}&\text{ if } c_n\to c\\
    1&\text{ if }c_n\to\infty
\end{cases}.
\] An even stronger result was given by Bollob\'as \cite{bollobas1984evolution}. He showed that for the random graph process, the hitting time for Hamiltonicity
is exactly the same as the hitting time for having minimum degree 2, that is, whp the very
edge which increases the minimum degree to 2 also makes the graph Hamiltonian.

In this paper, we take advantage of the study of the local resilience in random graphs and random digraphs, which was first introduced by Sudakov and Vu \cite{sudakov2008local}. 
Roughly speaking, the local resilience is the largest proportion of edges that one can delete from every vertex in a given graph $G$ satisfying a property $\mathcal P$, such that the resulting (sub)graph still satisfies $\mathcal P$. In \cite{sudakov2008local}, Sudakov and Vu showed that for any $\var > 0$, if $p> C\log n/n$ for some absolutely large constant $C>0$, then $G(n, p)$ typically has the property that every spanning subgraph with minimum degree at least $(1 + \var)np/2$ contains a perfect matching. Their result was then improved by Nenadov, Steger and Truji\'c \cite{nenadov2019resilience}, who proved the resilience result on random graph processes.

Besides the resilience on perfect matchings, Sudakov and Vu \cite{sudakov2008local} also study the resilience on Hamilton cycles. They showed that for any $\var > 0$, if $p$ is somewhat greater than $\log ^4n/n$, then $G(n, p)$ typically has the property that every spanning subgraph
with minimum degree at least $(1 + \var)np/2$ contains a Hamilton cycle. They also conjectured that this remains true as long as $p=(\log n+\omega(1))/n$, which was solved by Lee and Sudakov \cite{lee2012dirac}. Later, an even stronger result, the so-called ``hitting-time” statement, was shown by Nenadov, Steger and Truji\'c \cite{nenadov2019resilience}, and Montgomery \cite{montgomery2019hamiltonicity}, independently.

We shall use resilience results on perfect matchings in random bipartite graphs due to Sudakov and Vu~\cite{sudakov2008local} and Hamiltonicity in random digraphs by Montgomery~\cite{montgomery2020hamiltonicity}, see Lemma \ref{lem:bip-matching} and Lemma \ref{thm:resilience}, respectively.


\subsection{A rainbow setting}

In recent years, rainbow structures in graph systems have received a lot of attention. For instance, one might look for rainbow structures for matchings \cite{kupavskii2023rainbow,LU,lu2023better,  Montgomery22}, cycles \cite{aharoni2018rainbow,bradshaw2020transversals, bradshaw2022one,Cheng2,joos2020rainbow}, and pancyclicity \cite{bradshaw2020transversals,cheng2019rainbow}. There are also several attempts to generalize some traditional theorems into their rainbow settings, e.g.  ~\cite{ Howard2017A,  cheng2023rainbow}. Roughly speaking, to define a rainbow (hyper)graph, we need to start from a family of (hyper)graphs $\{G_1,\dots, G_m\}$ defined on the same vertex set, where each $c\in[m]$ is referred as a \emph{color}. A copy of a (hyper)graph $H$ with at most $m$ edges is called \emph{rainbow} if all the edges of $H$ are selected among the graphs $G_1,\dots, G_m$, and all the edges of $H$ have different colors. We define the rainbow graphs formally as below.

\begin{defi}\label{def:rainbow}
Let $\{G_1,\ldots,G_m\}$ be a family of not necessarily distinct (hyper)graphs on the same vertex set $V$. We say that a (hyper)graph $H$ on $V$ is \emph{rainbow} if there exists an injection $\varphi:E(H)\to[m]$ such that $e\in E(G_{\varphi(e)})$ for every $e\in E(H)$.
\end{defi}

A rainbow version of the Dirac-type problems in systems of graphs was conjectured by Aharoni et al.~\cite{aharoni2018rainbow} 
as follows: let $\{G_1,\ldots, G_{n}\}$ be a system of graphs on the same vertex set $V=[n]$ with minimum degree $\delta(G_c)\ge n/2$ for each $c\in[n]$. Then, there exists a rainbow Hamilton cycle. Cheng, Wang and Zhao \cite{cheng2019rainbow} verified the conjecture asymptotically and Joos and Kim \cite{joos2020rainbow} proved the full conjecture. Very recently, Bradshaw, Halasz, and Stacho \cite{bradshaw2022one} strengthened the result by showing that the system of graphs actually admit exponentially many rainbow Hamilton cycles under the same assumptions. Bradshaw \cite{bradshaw2020transversals} generalized the Dirac-type result for Hamiltonicity of bipartite graphs by Moon and Moser \cite{moon1963hamiltonian} to the rainbow setting.

In this note, we first give the Dirac-type result for rainbow perfect matchings which reads as follows.

\begin{thm}
\label{thm:dirac-rainbow-random-PM}
Let $\var\in(0,1)$ and $p=\omega\left(\log n/n\right)$ where $n\in 2\mathbb{N}$ is sufficiently large.
Let $\{G_1,\ldots, G_{n/2}\}$ be a family of graphs on the same vertex set $V=[n]$, where the graphs $G_c$s are independent samples of $G(n,p)$. 
Then, whp we have that for every spanning subgraphs $H_c\subseteq G_c$, $1\leq c\leq n/2$, with $\delta(H_c)\geq (\frac{1}{2}+\var)np$, the family  $\{H_1,\dots,H_{n/2}\}$ admits a rainbow perfect matching.
\end{thm}

Our proof of Theorem \ref{thm:dirac-rainbow-random-PM} is short and roughly relies on the following main idea. First, as a simple corollary of Chernoff's inequalities, the induced bipartite subgraph of $G(n,p)$ has concentrated vertex degrees for almost all balanced bipartition $[n]=V_1\cup V_2$. For such a bipartition, we pick a random permutation $\pi$ on $V_1$ and then define an auxiliary bipartite graph $B_\pi$ (see Definition \ref{def:auxiliary-graph}), which has the following crucial property: a perfect matching in $B_\pi$ corresponds to a rainbow perfect matching in the original family. Noting that the resulting bipartite subgraph of $G(n,p)$ can be regarded as a sample of random bipartite subgraph $B(n/2,n/2,p)$, the resilience result (see Lemma \ref{lem:bip-matching}) can thus ensure the resilience on perfect matchings for $B_\pi$, which concludes the proof by the crucial property of $B_\pi$.


With a similar approach except replacing the auxiliary bipartite graph by an auxiliary digraph, we also give the Dirac-type result for rainbow Hamilton cycles in random graphs as below.

\begin{thm}
\label{thm:dirac-rainbow-Hamilton-cycle}
Let $\var\in(0,1)$ and $p=\omega\left(\log n/n\right)$ where $n\in\mathbb{N}$ is sufficiently large.
Let  $\{G_1,\ldots, G_{n}\}$ be a family of graphs on the same vertex set $V=[n]$, where the graphs $G_c$s are independent samples of $G(n,p)$.  
Then, whp we have that for every spanning subgraphs $H_c\subseteq G_c$, $1\leq c\leq n$, with $\delta(H_c)\geq (\frac{1}{2}+\var)np$, the family $\{H_1,\dots,H_{n}\}$ admits a rainbow Hamilton cycle.
\end{thm}



Let us comment on the bound for $p$ in our main theorems. Unlike finding perfect matchings or Hamilton cycles in a single sample of $G(n,p)$, we only need one edge in each independent sample of $G(n,p)$ in order to find rainbow perfect matchings or rainbow Hamilton cycles. Thus, our bound $p=\omega(\log n/n)$ is probably not tight. However, if we take all the graphs $G_c$ to be the same sample of $G(n,p)$ as below, then our bound on $p$ is essentially tight.

\begin{thm}\label{thm:same-sample}
Let $\var\in(0,1)$ and $p=\omega\left(\log n/n\right)$ where $n\in\mathbb{N}$ is sufficiently large.
Let $G$ be a sample of $G(n,p)$ on vertex set $V=[n]$.  
Then, whp we have that for every spanning subgraphs $H_c\subseteq G$, $1\leq c\leq n$, with $\delta(H_c)\geq (\frac{1}{2}+\var)np$, the family $\{H_1,\dots,H_{n}\}$ admits a rainbow Hamilton cycle.
\end{thm}


Using an essentially the same proof as Theorem \ref{thm:dirac-rainbow-Hamilton-cycle}, one can also prove the analougous statement for perfect matchings as well whose proof is omitted.

\begin{thm}\label{thm:same-sample-PM}
Let $\var\in(0,1)$ and $p=\omega\left(\log n/n\right)$ where $n\in 2\mathbb{N}$ is sufficiently large.
Let $G$ be a sample of $G(n,p)$ on vertex set $V=[n]$. 
Then, whp we have that for every spanning subgraphs $H_c\subseteq G$, $1\leq c\leq n/2$, with $\delta(H_c)\geq (\frac{1}{2}+\var)np$, the family  $\{H_1,\dots,H_{n/2}\}$ admits a rainbow perfect matching.
\end{thm}

\subsection{Notation}
Given a graph $G$ and a subset 
$X\subseteq V(G)$, let $N(X)=\bigcup_{x\in X}N(x)$. 
For
two subsets $X,Y\subseteq V(G)$, we define $E_G(X,Y)$ to be the set of all edges $xy\in E(G)$ with $x\in X$ and $y\in Y$, and set $e_G(X,Y):=|E_G(X,Y)|$ (the subscript $G$ will be omitted whenever there is no risk of confusion). Moreover, $G[X,Y]$ is defined by the induced bipartite subgraph of $G$ with parts $X$ and $Y$.  When
$x \in V (G)$, $d_G(x)\coloneqq |N(x)|$ is the \emph{degree} of $x$ in $G$.
For a vertex $x\in V(G)$, and  a subset $Y\subset V(G)$, we define 
$d_G(x,Y)=|\{xy\in E(G):y\in Y\}|$. 
In particular, if $Y=\{y\}$ for some vertex $y\in V(G)$, then we write $d(x,y):=d(x,\{y\})$. For a graph $G$, we denote by $\delta(G)$ as its minimum degree, and $\Delta(G)$ as its maximum degree. For a digraph $D$, we denote $\delta^+(D),\delta^-(D)$ as its minimum out-degree and minimum in-degree, respectively. 
Moreover, let
\[
\delta^0(D)=\min\{\delta^+(D),\delta^-(D)\}.
\]

We will use various models of random (di)graphs defined as following. The 
random graph $G(n, p)$
on vertex set $[n] = \{1, \dots , n\}$ is a graph with all possible  edges chosen independently at random with probability $p$. Let $V_1$ and $V_2$ be disjoint vertex sets each of size $n/2$, where $n\in 2\mathbb N$.
The random bipartite graph $B(n/2,n/2,p)$ with parts $V_1$ and $V_2$  is a bipartite graph, such that for each $(i,j)\in V_1\times V_2$, $ij\in E(B(n/2,n/2,p))$ with probability $p$ and all edges $ij$ are chosen independently of each other.  The 
random digraph $D(n, p)$
on vertex set $[n]$ has all possible directed edges $(i,j)$, which are ordered pairs of vertices, chosen independently at random with probability $p$.

We write $a\in b\pm c$ as a shorthand for the double-sided inequality $b-c\leq a\leq b+c$, and write $a\notin b\pm c$ otherwise. If $f(n)/g(n) \to 0$ as $n \to\infty$, then we say $g(n) = \omega(f(n))$ and $f(n) = o(g(n))$. If there exists a
constant $C>0$ for which $|f(n)| \leq Cg(n)$ for all $n$, then we say $f(n) = O(g(n))$ and $g(n) = \Omega(f(n))$. If
$f = O(g(n))$ and $f(n) = \Omega(g(n))$, then we say that $f(n) = \Theta(g(n))$.




\section{Preliminary Results}\label{s:preliminary}

In this section, we collect some probabilistic tools and summarize some properties that typical graphs satisfy, which will be crucial in the proofs of our main results.

\subsection{Probabilistic tools}

We will use the following well-known bounds on the upper and lower tails of the binomial distribution, which is given by Chernoff (see Appendix A in \cite{AS}).

\begin{lem}[Chernoff's inequality for binomial distribution]\label{lem:Chernoff}
	Let $X$ be any sum of independent (not necessarily identical) Bernoulli random variables.
	\begin{enumerate}[label=(\roman*)]
		\item For $0 \leq \delta \leq 1$, we have
	\[	\Pr\left[X \leq (1 - \delta) \mathbb E[X]\right] \leq \exp\left\{-\frac{\delta^2 \mathbb E[X]}{2}\right\}.
    \]
		\item\label{right tail} For $\delta \geq 0$, we have
	\[	\Pr\left[X \geq (1 + \delta) \mathbb E[X]\right] \leq \exp\left\{-\frac{\delta^2 \mathbb E[X]}{2 + \delta}\right\}.\]
	\end{enumerate}
\end{lem}

The similar bounds can be also applied to hypergeometric distribution. More specifically, Hoeffding \cite{hoeffding1994probability} proved the following bounds (also see Section 23.5 in \cite{frieze2016introduction}).

\begin{lem}[Chernoff's inequality for hypergeometric distribution]\label{lem:hypergeometric}
Let 
$X\sim \mathrm{Hypergeometric}(N,K,n)$. Then
\begin{enumerate}[label=(\roman*)]
		\item For $0 \leq \delta \leq 1$, we have
	\[	\Pr\left[X \leq (1 - \delta) \mathbb E[X]\right] \leq \exp\left\{-\frac{\delta^2 \mathbb E[X]}{2}\right\}.
    \]
    \item\label{right tail hyper} For $\delta \geq 0$, we have
	\[	\Pr\left[X \geq (1 + \delta) \mathbb E[X]\right] \leq \exp\left\{-\frac{\delta^2 \mathbb E[X]}{2 + \delta}\right\}.\]
\end{enumerate}
\end{lem}

Our main probabilistic tool is the following concentration inequality, which is a special case of a more general result of McDiarmid \cite{mcdiarmid2002concentration} (see also Talagrand \cite{talagrand1995concentration}).

\begin{thm}\label{thm:Talagrand}
Let $c,r>0$ be constants, and let $h:S_n\to\mathbb{R}_{\geq0}$ be a function from the set of permutations on $[n]$ to nonnegative real numbers satisfying the following conditions:

\begin{enumerate}
    \item if $\pi,\pi'$ differ only in two places, then $|h(\pi)-h(\pi')|\leq c$;
    \item if $h(\pi) = s$, then there exists a subset $W\subseteq [n]$ of size at most $rs$, such that any $\pi'\in S_n$ that coincides with $\pi$ on $W$ satisfies $h(\pi')\geq s$. 
\end{enumerate}
Let $\pi\in S_n$ be a permutation chosen uniformly at random, and let $M\coloneqq M(h(\pi))$ be the median of the random variable $h(\pi)$. Then for each $t\geq 0$, we have
\[      
\Pr[h(\pi) \leq M - t] \leq 2 \exp{\left\{-\frac{t^2}{16rc^2M} \right\}}.
\]
\end{thm}

\subsection{Typical properties of graphs}

In this section, we collect some useful properties of a family of graphs on the same vertex set $V=[n]$. First, we show that for a family of $\mathrm{poly}(n)$ independent samples of random graph $G(n,p)$ on $[n]$, all the degrees in each graph are concentrated.

\begin{lem}\label{lem:similar-degrees}
Let $\var\in(0,1)$ and let $N=f(n)$ for some polynomial $f\geq 0$. Let $G_1,\ldots, G_{N}$ be independent samples of $G(n,p)$ on the same vertex set $V=[n]$. Then, whp we have 
\[
(1-\var)np\leq \delta(G_c)\leq \Delta(G_c) \leq (1+\var)np
\]
holds for all $c\in [N]$, provided that $p=\omega\left(\log n/n\right)$.
\end{lem}
\begin{proof}

Fix a vertex $v\in[n]$ and a color $c\in[N]$. Observe that $d_{G_c}(v) \sim \mathrm{Bin}(n-1, p)$, and therefore $\mu := \mathbb{E}[d_{G_c}(v)] = (n -1)p$.
Thus, since $p=\omega(\log n/n)$, by Lemma \ref{lem:Chernoff} we obtain that
\[
\Pr[d_{G_c}(v) \notin(1 \pm \var)\mu] \leq 2\exp{\left\{-\frac{\var^2\mu}{2+\var}\right\}}=o\left(\frac{1}{nN}\right).
\]
Taking a union bound over all vertices $v\in[n]$ and all colors $c\in[N]$, we conclude that
\[
\Pr[\exists v \in[n] ,\exists c\in[N]\text{ s.t. }d_{G_c}(v)\notin  (1\pm\var)\mu]=o(1).
\]
This completes the proof.
\end{proof}

Next, we show that given $n/2$ graphs $H_1,\dots, H_{n/2}$, if we take a random balanced bipartition of $[n]$, then whp all the induced bipartite subgraphs of $H_i$ on the bipartition have the ``correct'' degrees. 

\begin{lem}\label{lem:partition}
For every $\var>0$, there exists $C:=C(\var)>0$ for which the following holds. Let $n\in 2\mathbb{N}$ be sufficiently large, and let $m=n/2$. Let $H_1,\dots,H_m$ be graphs on the same vertex set $V=[n]$. Suppose that $\delta(H_c)\geq C\log_{}n$ for all $c\in [m]$. Then, $(1-o(1))$-fraction of balanced bipartitions $V=V_1\cup V_2$ satisfy the following property: for every vertex $v\in V$ and color $c\in[m]$, and for $i=1,2$ we have 
\[
d_{H_c}(v,V_i)\in(1\pm\var)\cdot \frac{d_{H_c}(v)}{2}.
\]
\end{lem}

\begin{proof}
Consider a uniformly random bipartition $V=V_1\cup V_2$ into sets both of size $m$. Fixed  a vertex $v\in[n]$ and a color $c\in[m]$,  note that $d_{H_c}(v,V_i)$ is hypergeometrically distributed with expected value $\frac{d_{H_c}(v)}{2}\pm1$ for each $i=1,2$. Therefore, by Lemma \ref{lem:hypergeometric} we obtain that
\[
\Pr\left[d_{H_c}(v,V_i)\notin (1\pm \var)\cdot \frac{d_{H_c}(v)}{2}\right]\leq 2\exp\left\{-\frac{\var^2 d_{H_c}(v)/2}{3}\right\}\leq 2\exp\left\{-\frac{C\var^2\log_{}n}{6}\right\}\leq 2n^{-3},
\]
where the last inequality holds for a large enough $C$.

Taking a union bound over all possible indices $i=1,2$, all vertices $v\in [n]$, and all colors $c\in[m]$, we conclude that
\[
\Pr\left[\exists i=1,2, \exists v \in[n] ,\exists c\in[m]\text{ s.t. }d_{H_c}(v,V_i)\notin (1\pm \var)\cdot \frac{d_{H_c}(v)}{2}\right]\leq 2nm\cdot 2n^{-3}=o(1).
\]
%
This completes the proof.
\end{proof}

\section{Resilience on rainbow perfect matchings}\label{s:matchings}

In this section, we prove Theorem \ref{thm:dirac-rainbow-random-PM}, which states the resilience result of rainbow perfect matchings in a family of independent samples of random graph $G(n,p)$. We will first construct an auxiliary bipartite graph with a crucial property: a perfect matching of the auxiliary bipartite graph corresponds to a 
rainbow perfect matching in the original family $\{H_1,\ldots,H_{n/2}\}$. We then study the resilience of the perfect matchings in this auxiliary bipartite graph, which further implies the corresponding resilience result on rainbow setting based on such property.

\subsection{An auxiliary bipartite graph}

Our proof uses an auxiliary bipartite graph as follows.

\begin{defi}\label{def:auxiliary-graph}
Let $n\in2\mathbb{N}$. Let $H_1',\dots,H_{n/2}'$ be bipartite graphs on the same vertex set $V=[n]$, each of which has the same bipartition $V=V_1 \cup V_2$ with $|V_1|=|V_2|=n/2$. By relabeling the vertices (if necessary), we may assume $V_1=[{n}/{2}]$. Given a permutation $\pi:V_1\to V_1$, the auxiliary bipartite graph $B_\pi:=B_\pi(H_1',\dots,H_{n/2}')$ is constructed as follows:
the parts of $B_{\pi}$ are $V_1$ and $V_2$; the edge set consists of all pairs $(i,j)\in V_1\times V_2$ such that $ij\in E(H_{\pi(i)}')$.
\end{defi}

\begin{rem}\label{rem:PM-implies-RM}
Observe that a perfect matching in $B_\pi$ corresponds to a rainbow perfect matching in the family $\{H_1', \dots, H_{n/2}'\}$. Indeed, every edge $ij$ in $B_\pi$ with $i\in V_1$ and $j\in V_2$ corresponds to an edge $ij$ in $H_{\pi(i)}'$, and since $\pi$ is a permutation of $V_1$, a perfect matching of $B_{\pi}$ uses exactly one edge from each $H'_i$.
\end{rem}

We also need the following result of Sudakov and Vu~\cite{sudakov2008local} on local resilience of perfect matchings on random bipartite graphs, whose proof can be found in the proof of~\cite[Theorem 3.1]{sudakov2008local}. Recall that random bipartite graph $B(n/2,n/2,p)$ defined on the partition $V_1\cup V_2$ with $|V_1|=|V_2|=n/2$ is a bipartite graph such that given $(i,j)\in V_1\times V_2$, $ij\in E(B(n/2,n/2,p))$ with probability $p$ and all possible edges $ij$ are chosen independently of each other.

\begin{lem}[Sudakov and Vu \cite{sudakov2008local}]
\label{lem:bip-matching}
Let $\var>0$ and $p=\omega\left(\log n/n\right)$ where $n\in 2\mathbb{N}$ is sufficiently large.
Let $V_1$ and $V_2$ be disjoint vertex sets each of size $n/2$.
Then, whp a spanning subgraph $H\subseteq B(n/2,n/2,p)$ defined on $V_1\cup V_2$ such that $\delta(H)\geq (\frac{1}{2}+\var)\frac{np}{2}$ contains a perfect matching.
\end{lem}

With the resilience result of perfect matchings on $B(n/2,n/2,p)$, we use the auxiliary bipartite graph to obtain the following bipartite version of Theorem~\ref{thm:dirac-rainbow-random-PM}.

\begin{thm}
\label{thm:dirac-rainbow-random-PM2}
Let $\var>0$ and $p=\omega\left(\log n/n\right)$ where $n\in 2\mathbb{N}$ is sufficiently large. Let $V_1$ and $V_2$ be disjoint vertex sets each of size $n/2$, and suppose that $G_1,\ldots, G_{n/2}$ are independent samples of $B(n/2,n/2,p)$ defined on $V_1\cup V_2$. 
Then, whp we have that for every spanning (bipartite) subgraphs $H_c\subseteq G_c$, $1\leq c\leq n/2$, with $\delta(H_c)\geq (\frac{1}{2}+\var)\frac{np}2$, the family $\{H_1,\dots,H_{n/2}\}$ admits a rainbow perfect matching.
\end{thm}

Before we prove Theorem \ref{thm:dirac-rainbow-random-PM2}, we first show Theorem \ref{thm:dirac-rainbow-random-PM} from the above theorem. To derive Theorem~\ref{thm:dirac-rainbow-random-PM}, it suffices to show that there is a balanced bipartition of the vertex set $[n]$ such that the bipartite subgraphs of our graphs inherit the degree condition, which was proved in Lemma~\ref{lem:partition}.

\begin{proof}[Proof of Theorem~\ref{thm:dirac-rainbow-random-PM}]
Let $\var > 0$ and $p \geq C\log n/n$, for a sufficiently large $C$. Let $m=n/2$. Let $G_1,\ldots, G_{m}$ be independent samples of $G(n,p)$ (on the same vertex set $V=[n]$).  
We wish to show that whp, for any subgraphs $H_c\subseteq G_c$ with $\delta(H_c) \geq (\frac{1}{2}+\var)np$, the family of graphs $\{H_1,\dots,H_m\}$ admits a rainbow perfect matching.

Recall that by Lemma \ref{lem:similar-degrees}, whp the family of graphs $G_1,\ldots, G_{m}$ satisfies
\begin{enumerate}[label=$(\dagger)$]
\item $(1-\var)np\leq \delta(G_c)\leq \Delta(G_c) \leq (1+\var)np$ for all $c\in[m]$. \label{item:dagger}
\end{enumerate}
For the rest of the proof, we condition on~\ref{item:dagger}.

Next, by Lemma \ref{lem:partition} with $\alpha \coloneqq\frac{\var}{100}$ 
in place of $\var$, we obtain that $1-o(1)$ fraction of balanced bipartitions $[n] = V_1\cup V_2$ satisfy the following: for every $v\in V_i$, $i=1,2$, and for every $c\in[m]$, we have 
\[
d_{H_c}(v,V_{3-i})\ge (1-\alpha)\cdot \frac{d_{H_c}(v)}{2} \ge (1-\alpha)\cdot \left(\frac{1}{2}+\var\right)\frac{np}2 \ge \left(\frac{1}{2}+\frac\var 2 \right)\frac{np}2.
\]
For such a partition, the host graph $G_c'\coloneqq G[V_1,V_2]$ of the subgraph $H_c$ is indeed a sample of $B(n/2,n/2,p)$, since the edges in the parts $V_1$ or $V_2$ are unrelated, while the edges between $V_1$ and $V_2$ have the same distribution as $G(n,p)$. Now, let $H_c':=H_c[V_1, V_2]$ be the spanning bipartite subgraphs of $H_c$ induced by the bipartition $V_1\cup V_2$.



Thus, for $1-o(1)$ fraction of balanced bipartitions, the graphs $H_c'\subseteq G_c'$ satisfy $\delta(H_c')\ge (\frac{1}{2}+\frac{\var}{2})\frac{np}{2}$, where $G_c'\sim B(n/2,n/2,p)$ are independent samples.
Therefore, by Theorem~\ref{thm:dirac-rainbow-random-PM2}, we obtain that by taking a random balanced bipartition $[n]=V_1\cup V_2$, whp the family $\{H_1',\dots,H_m'\}$ admits a rainbow perfect matching. Since $H_c'$ is an induced subgraph of $H_c$, the family $\{H_1,\ldots,H_m\}$ also admits a rainbow perfect matching. This completes the proof.
\end{proof}

\subsection{Most $B_{\pi}$'s have large minimum degree}

We are now left to prove Theorem \ref{thm:dirac-rainbow-random-PM2}. First, we prove that given bipartite graphs $H_1',\ldots,H'_m$ with a common balanced bipartition $V_1\cup V_2$ (each of size $m=n/2$) and large minimum degrees, the resulting auxiliary graph $B_\pi\coloneqq B_\pi(H_1',\dots,H_{m}')$ also has large minimum degree whp, where $\pi$ is a uniformly random permutation. The proof is a special case of Lemma 13 in \cite{FH}, but we include it for completeness. 

\begin{lem}\label{lem:large-degrees}
For every $0<\var<\frac{1}{2}$, there exists $C:=C(\var)$ for which the following holds for sufficiently large $n\in2\mathbb{N}$ and $p=C\log m/m$, where $m=n/2$. Let $H_1',\dots,H_m'$ be bipartite graphs on the same vertex set $V=[n]$ with the same balanced bipartition $V=V_1 \cup V_2$, where $V_1=[m]$. Suppose that $\delta(H_c')\geq (\frac{1}{2}+\var)mp$ for every $c\in[m]$. Let $\pi:[m]\to [m]$ be a uniformly random permutation, and let $B_\pi\coloneqq B_\pi(H_1',\dots,H_{m}')$. Then, whp we have $\delta(B_\pi)\geq (\frac{1}{2}+\frac{\var}{2})mp$.
\end{lem}

In order to use the McDiarmid's inequality (Theorem \ref{thm:Talagrand}) to prove the above lemma, we need to first bound the median of the random variable $d_{B_\pi}(j)$.

\begin{lem}\label{lem:median}
Let $0<\alpha<\frac{1}{2}$ and let $n\in2\mathbb{N}$ be sufficiently large. Let $m=n/2$. Let $H_1',\dots,H_{m}'$ be bipartite graphs on the same vertex set $V=[n]$ with the same balanced bipartitions $V=V_1 \cup V_2$, where $V_1=[m]$. Suppose that $\delta(H_c')\geq \frac{200}{\alpha^2}$ for all $c\in [m]$. 
Let $\pi$ be a uniformly random permutation on $V_1$. Let $B_\pi\coloneqq B_\pi(H_1',\dots,H_{m}')$, and let $\mu_i=\mathbb{E}[d_{B_\pi}(i)]$ for any $i\in V$. 
Then, for every $j\in V_2$, we have the median 
\[
M_j:=M(d_{B_\pi}(j))\in (1\pm\alpha)\mu_j.
\]
\end{lem}

\begin{rem}
The above lemma allows us to use $\mu_j$ instead of $M_j$ in Theorem \ref{thm:Talagrand} when it is applied to $d_{B_\pi}(j)$.
\end{rem}

\begin{proof}
Consider $B_\pi$, where $\pi$ is a uniformly random permutation on $V_1$. Let $j$ be some vertex in $V_2$. Let $\mu_j:=\mathbb{E}[d_{B_\pi}(j)]$ and $\sigma_j^2:= \mathrm{Var}(d_{B_\pi}(j))$. Moreover, for each $i\in V_1$, we define an indicator random variable $\mathds{1}_i$, where $\mathds{1}_i=1$ if $ij\in E(H_{\pi(i)}')$. Observe that $d_{B_{\pi}}(j)=\sum_{i=1}^m \mathds{1}_i.$ 

Applying Chebyshev's inequality, we have 
\[
\Pr[|d_{B_\pi}(j)-\mu_j|\geq \alpha\mu_j]\leq \frac{\sigma_j^2}{\alpha^2\mu_j^2}.
\]
If we can show that $\sigma_j^2\leq \frac{\alpha^2\mu_j^2}{100}$, then the result follows. Indeed, with probability at least 99/100 we have that $d_{B_\pi}(j)\in (1\pm\alpha)\mu_j$, and thus we conclude that the median $M_j$ also lies in this interval. Now the remaining part is to prove the desired inequality by computing $\mu_j=\mathbb{E}[d_{B_\pi}(j)]$ and $\sigma_j^2= \mathrm{Var}(d_{B_\pi}(j))$.

Note that the event $\left(\mathds{1}_i=1\right)$ only depends on the value of $\pi(i)$. There are $m$ possible values in total for $\pi(i)$, and exactly the colors in which $ij$ is an edge contribute to $\mathds{1}_i$. Recall that 
$d_{H_c'}(i,j)=1$ if $ij\in E(H_c')$, and $d_{H_c'}(i,j)=0$ otherwise. So
\[
\Pr[\mathds{1}_i=1]=\frac{\sum_{c=1}^m d_{H_c'}(i,j)}{m}.
\] 
By linearity of expectations, we have
\[
\mu_j=\sum_{i=1}^m \mathbb{E}[\mathds{1}_i]=\sum_{i=1}^m\frac{\sum_{c=1}^m d_{H_c'}(i,j)}{m} =\sum_{c=1}^m\frac{\sum_{i=1}^m d_{H_c'}(i,j)}{m}=\sum_{c=1}
^m \frac{d_{H_c'}(j)}{m}.
\]
To compute the variance, note that for each $i\neq k$ in $V_1$, we have
\[
\begin{aligned}
\mathbb{E}\left[\mathds{1}_i\mathds{1}_k\right]
&=\sum_{c=1}^m \Pr\left[\mathds{1}_i=\mathds{1}_k=1|\pi(i)=c\right] \Pr\left[\pi(i)=c\right]\\
&=\sum_{c=1}^m \frac{1}{m} d_{H_c'}(i,j)\Pr\left[\mathds{1}_k=1|\pi(i)=c\right]\\
&=\sum_{c=1}^m \frac{1}{m} d_{H_c'}(i,j)\frac{\sum_{c'\neq c} d_{H_{c'}'}(k,j)}{m-1}\\
&\leq \frac{m}{m-1}\sum_{c=1}^m \frac{d_{H_c'}(i,j)}{m}\sum_{c'=1}^m\frac{ d_{H_{c'}'}(k,j)}{m}\\
&=\frac{m}{m-1}\mathbb{E}[\mathds{1}_j]\mathbb{E}[\mathds{1}_k].
\end{aligned}
\]
Thus,
\[
\begin{aligned}
\mathrm{Var}\left(d_{B_\pi}\left(j\right)\right)&=\mathrm{Var}\left(\sum_{i=1}^m \mathds{1}_i\right)
=\sum_{i=1}^m \mathrm{Var}\left(\mathds{1}_i\right) + \sum_{i\neq k}\mathrm{Cov}\left(\mathds{1}_i,\mathds{1}_k\right)\\
&\leq \mu_j+\sum_{i\neq k}\left(\mathbb{E}\left[\mathds{1}_i\mathds{1}_k\right] - \mathbb{E}\left[\mathds{1}_i\right]\mathbb{E}\left[\mathds{1}_k\right]\right)\\
&\leq \mu_j+\frac{1}{m-1}\sum_{i\neq k}\mathbb{E}[\mathds{1}_i]\mathbb{E}[\mathds{1}_k]\\
&\leq \mu_j+\frac{1}{m-1}\sum_{i=1}^m\mathbb{E}[\mathds{1}_i]\sum_{k=1}^m\mathbb{E}[\mathds{1}_k]\\
&=\mu_j+\frac{1}{m-1}\mu_j^2.
\end{aligned}
\]
To complete the proof, first observe that we have $\frac{1}{m-1}\mu_j^2\leq \frac{\alpha^2\mu_j^2}{200}$ since $m$ is sufficiently large. Also, we have $\mu_j\leq \frac{\alpha^2\mu_j^2}{200}$ since $\mu_j\geq \frac{200}{\alpha^2}$ by assumption. 
Now we obtain $\sigma_j^2\leq \frac{\alpha^2\mu_j^2}{100}$ and the lemma follows.
\end{proof}

Now we are able to prove Lemma \ref{lem:large-degrees}.

\begin{proof}[Proof of Lemma \ref{lem:large-degrees}]
Consider $B_\pi$, where $\pi$ is a uniformly random permutation. As $\delta(H_c')\geq (\frac{1}{2}+\var)mp$ for every $c\in[m]$, it is guaranteed that for all $i\in V_1$ we have that $d_{B_\pi}(i)\geq (\frac{1}{2}+\var)mp$.
Now fix some $j\in V_2$ and observe from the proof of Lemma \ref{lem:median}, under the same notation, that $\mu_j:=\mathbb{E}[d_{B_\pi}(j)]\geq (\frac{1}{2}+\var)mp.$

In order to complete the proof, we want to show that the $d_{B_\pi}(j)$s are highly concentrated using Theorem \ref{thm:Talagrand}. To this end, let $h(\pi):=d_{B_\pi}(j)$ and note that swapping any two elements of $\pi$ can change the value of $h$ by at most 2. Moreover, note that if $h(\pi) =d_{B_\pi}(j)= s$, then it is enough to choose the subset $W\coloneqq N_{B_\pi}(j)$ so that every permutation $\pi'$ that coincides $\pi$ on $W$ satisfies $h(\pi')\geq s$. Therefore, $h(\pi)$ satisfies the conditions of Theorem \ref{thm:Talagrand} with $c = 2$ and $r = 1$.

Now, let $\alpha=\frac{\var}{100}$, and observe that by Lemma \ref{lem:median} the median $M$ of $d_{B_\pi}(j)$ lies in the interval $(1\pm\alpha)\mu_j$.
Therefore, we have
\[
\Pr \left[h(\pi)\leq \left(\frac{1}{2}+\frac{\var}{2}\right)mp\right]\leq \Pr\left[h(\pi)\leq \left(1-\frac{\var}{2}\right)\mu_j\right]
\]
and the latter can be further upper bounded by
\[
\Pr\left[h(\pi)\leq \left(1-\frac{\var}{2}\right)(1-\alpha)^{-1}M\right]\leq \Pr\left[h(\pi)\leq \left(1-\frac{\var}{4}\right)M\right].
\]
Now, by Theorem \ref{thm:Talagrand} we obtain that
\[
\Pr\left[h(\pi)\leq \left(\frac{1}{2}+\frac{\var}{2}\right)mp\right]\leq 2\exp{\left\{-\frac{(\var M/4)^2}{64M}\right\}}.
\]
Next, using (again) the fact that $M\in (1\pm\alpha)\mu_j$ and that $\mu_j=\Theta(mp)\geq C\log_{}m$, we can upper bound the above right hand side by
$2\exp{ (-\Theta(mp))} \leq n^{-2}$.
Finally, in order to complete the proof, we take a union bound over all $j\in V_2$ and conclude that whp $\delta(B_\pi)\geq (\frac{1}{2}+\frac{\var}{2})mp$.
\end{proof}

\subsection{Proof of Theorem \ref{thm:dirac-rainbow-random-PM2}}


Now we are ready to prove Theorem~\ref{thm:dirac-rainbow-random-PM2}.
\begin{proof}[Proof of Theorem~\ref{thm:dirac-rainbow-random-PM2}]
Let $\var > 0$ and $p \geq C\log n/n$, for a sufficiently large $C$. Let $m=n/2$. Let $G_1,\ldots, G_{m}$ be independent samples of $B(n/2,n/2,p)$ on $V_1\cup V_2$.  
%
We wish to demonstrate that whp, for any (bipartite) subgraphs $H_c\subseteq G_c$ with $\delta(H_c) \geq (\frac{1}{2}+\var)mp$, the family of graphs $\{H_1,\dots,H_m\}$ has a rainbow perfect matching.

First, recall that Lemma \ref{lem:large-degrees} (with input graphs $H_1,\dots, H_m$) guarantees that for almost all permutation $\pi$ on $V_1$, the auxiliary bipartite graph $H_\pi:=B_\pi(H_1,\dots, H_m)$ (as defined in Definition \ref{def:auxiliary-graph}) satisfies 
\begin{enumerate}[label=$(\ddagger)$]
\item $\delta(H_\pi)\geq (\frac{1}{2}+\frac{\var}{2})mp$. 
\label{item:ddagger}
\end{enumerate}
We focus on all $\pi$ satisfying above and  condition on~\ref{item:ddagger}. 

Note that the bipartite graph $G_\pi:=B_\pi(G_1,\dots, G_m)$ is a random bipartite graph, where all pairs in $V_1\times V_2$ are presented with probability $p$, and independent with other pairs of vertices.
Moreover, by definition, $H_\pi$ is a subgraph of $G_\pi$.
Thus, by~\ref{item:ddagger}, we have that $H_\pi$ is a subgraph of $G_\pi$ with $\delta(H_\pi)\geq (\frac{1}{2}+\frac{\var}{2})mp$.
Therefore, by Lemma~\ref{lem:bip-matching}, whp $H_\pi$ contains a perfect matching, which by Remark \ref{rem:PM-implies-RM} implies that the family $\{H_1,\dots,H_m\}$ whp admits a rainbow perfect matching. This completes the proof.
%
%
\end{proof}

\section{Resilience on rainbow Hamilton cycles}\label{s:cycles}
In this section, we prove Theorem \ref{thm:dirac-rainbow-Hamilton-cycle}, which states the resilience result of rainbow Hamilton cycles in a family of independent samples of random graph $G(n,p)$. The proof is quite similar to the proof of Theorem \ref{thm:dirac-rainbow-random-PM}, where the only difference is the construction of an auxiliary digraph rather than a bipartite graph. In fact, we will construct the auxiliary digraph whose Hamilton cycle corresponds to a rainbow Hamilton cycle in the original family $\{H_1,\ldots,H_{n}\}$. Then we study the resilience of the Hamilton cycles in this auxiliary digraph. We will also give a proof sketch of Theorem \ref{thm:same-sample}, whose proof is quite similar to Theorem \ref{thm:dirac-rainbow-Hamilton-cycle}.

\subsection{An auxiliary digraph}
Similar to Section \ref{s:matchings}, we define an auxiliary digraph as follows.

\begin{defi}\label{def:auxiliary-digraph}
Let $H_1',\dots,H_{n}'$ be graphs on the same vertex set $V=[n]$. Given a permutation $\pi:V\to V$, the auxiliary digraph ${D}_\pi:={D}_\pi(H_1',\dots,H_{n}')$ is constructed as follows: \\
$V$ is the vertex set of ${D}_\pi$ and for any pair of vertices $(i,j)\in V\times V$, $(i,j)\in E({D}_\pi)$ if and only if $ij\in E(H_{\pi(i)}')$.
\end{defi}

\begin{rem}\label{rem:auxiliary-digraph}
Observe that a directed Hamilton cycle in the auxiliary digraph ${D}_\pi$ corresponds to a rainbow Hamilton cycle in the family $\{H_1',\dots,H_{n}'\}$. Indeed, a directed Hamilton cycle in ${D}_\pi$ is a directed Hamilton cycle in $K_n$ whose edge $(i,j)$ belongs to distinct $H_{\pi(i)}'$ since $\pi$ is a permutation.
\end{rem}

We also need the following result on local resilience of Hamiltonicity in random digraphs due to Montgomery \cite{montgomery2020hamiltonicity} (in fact, Montgomery proved a way stronger result but the following is enough for our needs). 


\begin{lem}[Montgomery \cite{montgomery2020hamiltonicity}]\label{thm:resilience}
Let $\var>0$. Then whp a spanning subdigraph $D\sim D(n,p)$ defined on $V=[n]$ such that $\delta^0(D)\geq (\frac{1}{2}+\var)np$ contains a Hamilton cycle, provided that $p=\omega\left(\log n/n\right)$.
\end{lem}

\subsection{Most $D_{\pi}$'s have large minimum degree}

The following lemma says that given digraphs $H_1',\ldots,H'_m$ with large minimum semidegrees (minimum of out-degrees and in-degrees), the resulting auxiliary digraph $D_\pi$ also has large minimum degree whp, where $\pi$ is a uniformly random permutation.

\begin{lem}
\label{lem:large-degrees-digraphs}
For every $\var>0$, there exists $C:=C(\var)$ for which the following holds for sufficiently large $n\in\mathbb{N}$ and $p=C\log n/n$. Let $H_1',\dots,H_n'$ be graphs on the same vertex set $V=[n]$. Suppose that $\delta(H_c')\geq (\frac{1}{2}+\var)np$ for every $c\in[n]$. Let $\pi:V\to V$ be a uniformly random permutation, and let $D_\pi:=D_\pi(H_1',\dots, H_{n}')$. Then, whp we have $\delta^0(D_\pi)\geq (\frac{1}{2}+\frac{\var}{2})np$.
\end{lem}

The proof of Lemma \ref{lem:large-degrees-digraphs} is very similar to that of Lemma \ref{lem:large-degrees} so we leave it to the appendix.

\subsection{Proof of Theorem \ref{thm:dirac-rainbow-Hamilton-cycle}}

Now we give a proof of Theorem \ref{thm:dirac-rainbow-Hamilton-cycle}.


\begin{proof}[Proof of Theorem~\ref{thm:dirac-rainbow-Hamilton-cycle}]

Let $\var > 0$ and $p \geq C\log n/n$, for a sufficiently large $C$. Let $G_1,\ldots, G_{n}$ be independent samples of $G(n,p)$ (on the same vertex set $V=[n]$). We wish to show that whp, for any subgraphs $H_c\subseteq G_c$ with $\delta(H_c) \geq (1/2+\var)np$, the family of graphs $H_1,\dots,H_n$ has a rainbow Hamilton cycle.

Similar to the proof of Theorem \ref{thm:dirac-rainbow-random-PM}, whp the family of graphs $G_1,\ldots, G_{n}$ satisfies \ref{item:dagger} for $m=n$. For the rest of the proof, we condition on~\ref{item:dagger}.

Now, let $\pi$ be a permutation on $V$ chosen uniformly at random. Note that the digraph $G_\pi:=D_\pi(G_1,\dots, G_n)$ is a random digraph, where all pairs are presented with probability $p$, and independent with other pairs of vertices. By exposing the subgraphs $H_c\subseteq G_c$, Lemma \ref{lem:large-degrees-digraphs} (with input graphs $H_1,\dots, H_n$) guarantees that whp $H_\pi:=D_\pi(H_1,\dots, H_n)$ (as defined in Definition \ref{def:auxiliary-digraph}) satisfies 
\begin{enumerate}[label=$(*)$]
\item $\delta^0(H_\pi)\geq (\frac{1}{2}+\frac{\var}{2})np$. \label{item:dag}
\end{enumerate}

Therefore, by Lemma~\ref{thm:resilience}, whp $H_\pi$ contains a Hamilton cycle, which by Remark \ref{rem:auxiliary-digraph} implies that the family $\{H_1,\dots,H_n\}$ whp admits a rainbow Hamilton cycle. This completes the proof.
\end{proof}

\subsection{Proof of Theorem \ref{thm:same-sample}}

We only give a proof sketch of Theorem \ref{thm:same-sample}, since its proof is quite similar to the proof of Theorem \ref{thm:dirac-rainbow-Hamilton-cycle}.

\begin{proof}[Proof sketch of Theorem \ref{thm:same-sample}]
Let $\var > 0$ and $p \geq C\log n/n$, for a sufficiently large $C$. Let $G$ be a sample of $G(n,p)$ on vertex set $V=[n]$. Consider the auxiliary digraph $G_\pi\coloneqq D_\pi(G,\ldots,G)$. Noticing that all the directed edges in $G_\pi$ are in the same color, we have that $D_\pi$ is a random digraph, where all pairs are presented with probability $p$, and independent with other pairs of vertices. By exposing the subgraphs $H_c\subseteq G$, Lemma \ref{lem:large-degrees-digraphs} (with input graphs $H_1,\dots, H_n$) guarantees that whp $H_\pi:=D_\pi(H_1,\dots, H_n)$ (as defined in Definition \ref{def:auxiliary-digraph}) satisfies 
\begin{enumerate}[label=$(*)$]
\item $\delta^0(H_\pi)\geq (\frac{1}{2}+\frac{\var}{2})np$. \label{item:dag}
\end{enumerate}

Therefore, by Lemma~\ref{thm:resilience}, whp $H_\pi$ contains a Hamilton cycle, which by Remark \ref{rem:auxiliary-digraph} implies that the family $\{H_1,\dots,H_n\}$ whp admits a rainbow Hamilton cycle. This completes the proof. 
\end{proof}

\bibliographystyle{abbrv}
\bibliography{references}
\appendix
\section{Proof of Lemma \ref{lem:large-degrees-digraphs}}
In this section, we finish the proof of Lemma \ref{lem:large-degrees-digraphs} which we omitted in Section \ref{s:cycles}. To apply McDiarmid's inequality (Theorem \ref{thm:Talagrand}), we first give a bound on the median of the random variable $d_{D_\pi}^-(i)$.

\begin{lem}\label{lem:median-digraphs}
Let $0<\alpha<\frac{1}{2}$ and let $n\in\mathbb{N}$ be sufficiently large. Let $H_1',\dots,H_n'$ be graphs on the same vertex set $V=[n]$. Suppose that $\delta\left(H_c'\right)\geq \frac{200}{\alpha^2}$ for all $c\in[n]$, Let $\pi$ be a uniformly random permutation on $V$ and $\mu_i=\mathbb{E}\left[d_{D_\pi}^-\left(i\right)\right]$. Then, for every $i\in V$, we have
\[
M_i:=M(d^-_{D_\pi}(i))\in \left(1\pm\alpha\right)\mu_i.
\]
\end{lem}

\begin{rem}
The above lemma allows us to use $\mu_i$ instead of $M_i$ in Theorem \ref{thm:Talagrand} when it is applied to $d^-_{D_\pi}(i)$.
\end{rem}

\begin{proof}
Consider $D_\pi$, where $\pi$ is a uniformly random permutation on $V$. Let $i$ be some vertex in $V$. Let $\mu_i:=\mathbb{E}\left[d^-_{D_\pi}(i)\right]$ and $\sigma_i^2:= \mathrm{Var}(d^-_{D_\pi}(i))$. Moreover, for each $j\in V$, we define an indicator random variable $\mathds{1}_j$, where $\mathds{1}_j=1$ if $ij\in E(H_{\pi(j)}')$. Observe that $d^-_{D_\pi}(i)=\sum_{j=1}^n \mathds{1}_j.$ 

Applying Chebyshev's inequality, we have 
\[
\Pr[|d^-_{D_\pi}(i)-\mu_i|\geq \alpha\mu_i]\leq \frac{\sigma_i^2}{\alpha^2\mu_i^2}.
\]
If we can show that $\sigma_i^2\leq \frac{\alpha^2\mu_i^2}{100}$, then the result follows. Indeed, with probability at least 99/100 we have that $d^-_{D_\pi}(i)\in (1\pm\alpha)\mu_i$ and thus we conclude that the median $M_i$ also lies in this interval. Now the remaining part is to prove the desired inequality by computing $\mu_i=\mathbb{E}[d^-_{D_\pi}(i)]$ and $\sigma_i^2= \mathrm{Var}(d^-_{D_\pi}(i))$.

Note that the event $\left(\mathds{1}_j=1\right)$ only depends on the value of $\pi(j)$. There are $n$ possible values in total for $\pi(j)$, and exactly all of the colors in which $ij$ is an edge contributes to $\mathds{1}_j$. Recall that $d_{H_c'}(i,j)=1$ if $ij\in E(H_c')$, and $d_{H_c'}(i,j)=0$ otherwise. So
\[
\Pr[\mathds{1}_j=1]=\frac{\sum_{c=1}^n d_{H_c'}(i,j)}{n}.
\]
By linearity of expectations, we have
\[
\mu_i=\sum_{j=1}^{n} \mathbb{E}[\mathds{1}_j]=\sum_{j=1}^{n}\frac{\sum_{c=1}^n d_{H_c'}(i,j)}{n} =\sum_{c=1}^n\frac{\sum_{j=1}^{n} d_{H_c'}(i,j)}{n}=\sum_{c=1}
^n \frac{d_{H_c'}(i)}{n}.
\]
To compute the variance, note that for each $j\neq k$ in $V$, we have
\[
\begin{aligned}
\mathbb{E}[\mathds{1}_j\mathds{1}_k]
&=\sum_{c=1}^n \Pr[\mathds{1}_j=\mathds{1}_k=1|\pi(j)=c]\Pr[\pi(j)=c]\\
&=\sum_{c=1}^n \frac{1}{n} d_{H_c'}(i,j)\Pr[\mathds{1}_k=1|\pi(j)=c]\\
&=\sum_{c=1}^n \frac{1}{n} d_{H_c'}(i,j)\frac{\sum_{c'\neq c} d_{H_{c'}'}(i,k)}{n-1}\\
&\leq \frac{n}{n-1}\sum_{c=1}^n \frac{d_{H_c'}(i,j)}{n}\sum_{c'=1}^n\frac{ d_{H_{c'}'}(i,k)}{n}\\
&=\frac{n}{n-1}\mathbb{E}[\mathds{1}_j]\mathbb{E}[\mathds{1}_k].
\end{aligned}
\]
Thus,
\[
\begin{aligned}
\mathrm{Var}(d^-_{D_\pi}(i))&=\mathrm{Var}(\sum_{j=1}^n \mathds{1}_j)
=\sum_{j=1}^n \mathrm{Var}(\mathds{1}_j) + \sum_{j\neq k}\mathrm{Cov}(\mathds{1}_j,\mathds{1}_k)\\
&\leq \mu_i+\sum_{j\neq k}(\mathbb{E}[\mathds{1}_j\mathds{1}_k] - \mathbb{E}[\mathds{1}_j]\mathbb{E}[\mathds{1}_k])\\
&\leq \mu_i+\frac{1}{n-1}\sum_{j\neq k}\mathbb{E}[\mathds{1}_j]\mathbb{E}[\mathds{1}_k]\\
&\leq \mu_i+\frac{1}{n-1}\sum_{j=1}^n\mathbb{E}[\mathds{1}_j]\sum_{k=1}^n\mathbb{E}[\mathds{1}_k]\\
&=\mu_i+\frac{1}{n-1}\mu_i^2.
\end{aligned}
\]
To complete the proof, first observe that we have $\frac{1}{n-1}\mu_i^2\leq \frac{\alpha^2\mu_i^2}{200}$ since $n$ is sufficiently large. Also, we have $\mu_i\leq \frac{\alpha^2\mu_i^2}{200}$ since $\mu_i\geq \frac{200}{\alpha^2}$ by assumption. 
Now we obtain $\sigma_i^2\leq \frac{\alpha^2\mu_i^2}{100}$ and the lemma follows.
\end{proof}

Now we are ready to prove Lemma \ref{lem:large-degrees-digraphs}.

\begin{proof}[Proof of Lemma \ref{lem:large-degrees-digraphs}]
Consider $D_\pi$, where $\pi$ is a uniformly random permutation on $V=[n]$. As $\delta\left(H_c'\right)\geq \left(\frac{1}{2}+\var\right)np$ for every $c\in [n]$ by assumption, it is guaranteed that $\delta^+\left(D_\pi\right)\geq \left(\frac{1}{2}+\var\right)np$. So it suffices to prove that $\delta^-\left(D_\pi\right)\geq \left(\frac{1}{2}+\var\right)np$. Now fix some vertex $i\in V$ and observe from the proof of Lemma \ref{lem:median-digraphs}, under the same notations, that $\mu_i:=\mathbb{E}\left[d_{D_\pi}^-(i)\right]\geq \left(\frac{1}{2}+\var\right)np$. 

In order to complete the proof, we want to show that the $d^-_{D_\pi}(i)$s are ``highly concentrated" using Theorem \ref{thm:Talagrand}. To this end, let $h(\pi):=d^-_{D_\pi}(i)$ and note that swapping any two elements of $\pi$ can change the value of $h$ by at most 2. Moreover, note that if $h(\pi) =d^-_{D_\pi}(i)= s$, then it is enough to choose the subset $W:=N^-_{D_\pi}(i)$ so that every permutation $\pi'$ that coincides $\pi$ on $W$ satisfies $h(\pi')\geq s$. Therefore, $h(\pi)$ satisfies the conditions of Theorem \ref{thm:Talagrand} with $c = 2$ and $r = 1$.

Now, let $\alpha=\frac{\var}{100}$, and observe that by Lemma \ref{lem:median-digraphs} the median $M$ of $d^-_{D_\pi}(i)$ lies in the interval $(1\pm\alpha)\mu_i$.
Therefore, we have
\[
\Pr \left[h(\pi)\leq \left(\frac{1}{2}+\frac{\var}{2}\right)np\right]\leq \Pr\left[h(\pi)\leq \left(1-\frac{\var}{2}\right)\mu_i\right]
\]
and the latter can be further upper bounded by
\[
\Pr\left[h(\pi)\leq \left(1-\frac{\var}{2}\right)(1-\alpha)^{-1}M\right]\leq \Pr\left[h(\pi)\leq \left(1-\frac{\var}{4}\right)M\right].
\]
Now, by Theorem \ref{thm:Talagrand} we obtain that
\[
\Pr\left[h(\pi)\leq \left(\frac{1}{2}+\frac{\var}{2}\right)np\right]\leq 2\exp{\left\{-\frac{(\var M/4)^2}{64M}\right\}}.
\]
Next, using (again) the fact that $M\in (1\pm\alpha)\mu_i$ and that $\mu_i=\Theta(np)\geq C\log_{}n$, we can upper bound the above right hand side by
$2\exp{ (-\Theta(np))} \leq n^{-2}$.
Finally, in order to complete the proof, we take a union bound over all $i\in V$ and conclude that whp $\delta^-(D_\pi)\geq (\frac{1}{2}+\frac{\var}{2})np$.
\end{proof}
\end{document}